\documentclass[12pt]{article}
\usepackage{amsmath,amssymb}
\input{epsf}

\newtheorem{theorem}{Theorem}
\newtheorem{lemma}[theorem]{Lemma}
\newtheorem{claim}[theorem]{Claim}

\newtheorem{corollary}[theorem]{Corollary}

\newenvironment{proof}{{\bf Proof.}}{\hfill\rule{2mm}{2mm}}

\newtheorem{remarka}[theorem]{Remark}
\newenvironment{remark}{\begin{remarka}\rm}{\end{remarka}}

\newtheorem{prelem}{{\bf Theorem}}

\newenvironment{lem}{\begin{prelem}{\hspace{-0.5
               em}{\bf.}}}{\end{prelem}}

\newtheorem{preconj}[prelem]{{\bf Conjecture}}

\newenvironment{conj}{\begin{preconj}{\hspace{-0.5
               em}{\bf.}}}{\end{preconj}}

\def\Var {{\rm Var}}
\def\Re {{\rm Re}}

\title{\bf Fourier analysis and large independent sets in powers of complete graphs.}
\author{
{\bf Mahya  Ghandehari$^a$ and Hamed Hatami$^b$} \\
$^a${\small\it Department of Pure Mathematics}\\
{\small {University of Waterloo}}\\
$^b${\small\it Department of Computer Science }\\
{\small University of Toronto}}
\date{}
\begin{document}
\maketitle
\begin{abstract}
For constant $r$ and arbitrary $n$, it was known that in the graph
$K_r^n$ any independent set of size close to the maximum is close
to some independent set of maximum size. We prove that this
statement holds for arbitrary $r$ and $n$.
\end{abstract}
\newline
{{\sc Keywords:} independent sets, weak product, Fourier analysis,
Harmonic analysis.

\section{Introduction \label{intro}}
The \emph{weak product} of $G$ and $H$, denoted by $G\times H$ is
defined as follows: The vertex set of $G\times H$  is the
Cartesian product of the vertex sets of $G$ and $H$. Two vertices
$(g_1,h_1)$ and $(g_2,h_2)$ are adjacent in $G\times H$ if
$g_1g_2$ is an edge of $G$ and $h_1h_2$ is an edge of $H$.

In this paper we consider the product of complete graphs on $r>2$ vertices,
$$G=K_r^n=\times_{j=1}^n K_r.$$
We identify the vertices of $G$ with the elements of $\mathbb{Z}_r^n$. By the definition of product,
two vertices are adjacent in $G$ iff the corresponding vectors differ in every coordinate.

Let $0\leq i \leq r-1$ and $1\leq j\leq n$ be two fixed integers.
It is obvious that the set of all vertices of $G$ which has $i$ in
the $j$th coordinate forms an independent set. In fact, for $r>2$,
these sets are the only maximum independent sets of
$G$~\cite{lovazs}. A generalization of this result has been shown
in \cite{ADFS} through the following Theorem:
\begin{lem}
{\bf \cite{ADFS}} \label{t2} For every $r\geq 3$, there exists a
constant $M=M(r)$ such that for any $\epsilon>0$ the following is
true. Let $G=K^n_r$ and $J$ be an independent set such that
$\frac{|J|}{|G|}=\frac{1}{r}(1-\epsilon)$. Then there exists an
independent set $I$ with $\frac{|I|}{|G|}=\frac{1}{r}$ such that
$\frac{|J\triangle I|}{|G|}<\frac{M\epsilon}{r}$.
\end{lem}

In Theorem \ref{t2}, ``$\triangle$'' denotes the symmetric
difference. Theorem \ref{t2} asserts that any independent set
which is close to being of maximum size is close to being
determined by one coordinate. The function $M(r)$ that is obtained
in~\cite{ADFS} depends on $r$. When $r$ is a constant, for every
constant $\delta>0$ one can choose $\epsilon$ to be a sufficiently
small constant so that $\frac{|J\triangle
I|}{|G|}<\frac{\delta}{r}$. But when $r$ tends to infinity, to
obtain  any nontrivial result from Theorem~\ref{t2}, $\epsilon$
must be less than $\frac{1}{M(r)}$ which is not a constant. The
main result of this paper is to show that in Theorem~\ref{t2}, $M$
does not need to be a function of $r$. Note that this major
improvement makes Theorem~\ref{t2} as powerful for large values of
$r$ as for constant $r$. We formalize this in the following
theorem.

\begin{theorem}
\label{th1} Let $G=K^n_r$, $r\geq 20$ and $\epsilon<10^{-9}$.
Suppose that $J$ is an independent set of $G$ such that
$\frac{|J|}{|G|}=\frac{1}{r}(1-\epsilon)$. Then there exists an
independent set $I$ with $\frac{|I|}{|G|}=\frac{1}{r}$ such that
$\frac{|J\triangle I|}{|G|}<\frac{40\epsilon}{r}$.
\end{theorem}

\begin{remark}
\label{remark1} Note that for $\epsilon\ge 10^{-9}$, we have the
trivial bound $\frac{|J\triangle I|}{|G|}\le \frac{2 \times
10^9\epsilon}{r}$, where $I$ is an arbitrary independent set. We
also assumed that $r \ge 20$, for some technical reasons. However
one can use Theorem~\ref{t2} when $r<20$, as $M(r)$ is a constant
for those values of $r$.
\end{remark}

Let $I$ be a maximum independent set of $G=K^n_r$, and $J$ be an
independent set of $G$ such that $J \not\subseteq I$. Then
obviously, $\frac{|I\setminus J|}{|G|} \geq
\frac{(r-1)^{n-1}}{r^n}$. So we obtain the following as a
corollary of Theorem~\ref{th1}.

\begin{corollary}
\label{cor1} Let $G=K^n_r$, $r\geq 20$ and $\epsilon<c$ where
$c=\min(10^{-9},(1-\frac{1}{r})^{n-1})/40$. Let $J$ be an
independent set such that
$\frac{|J|}{|G|}=\frac{1}{r}(1-\epsilon)$. Then there exists an
independent set $I$ with $\frac{|I|}{|G|}=\frac{1}{r}$ such that
$J \subseteq I$.
\end{corollary}

Note that if in Corollary~\ref{cor1}, $r>c'n$ for some constant $c'$,
then one can take $c$ to be a constant that does not depend on $n$.

The proof of Theorem \ref{th1} as well as Theorem \ref{t2} is
based on Fourier analysis on the group $\mathbb{Z}_r^n$. Fourier
analysis has been shown to be very useful in the study of Boolean
functions. One can refer to \cite{ADFS, AKRS, ALM, BS, BKS, B1,
F,FKN, KKL, KS, LMN, Mesh, T} to see some examples. In order to
prove Theorem~\ref{th1} we show that a Boolean function which has
most of its 2-norm weight concentrated on the first two
levels\footnote{Defined formally below} of its Fourier expansion
is close to being determined by one coordinate. Thus
Lemma~\ref{l1} which formulates this might be of independent
interest as a result in the direction of extending results
of~\cite{FKN,B1,KS} from $\mathbb{Z}_2^n$ to $\mathbb{Z}_r^n$.

Section 2 is devoted to a very brief introduction to Fourier
analysis of $\mathbb{Z}_r^n$ and introducing notations and some of
the necessary tools. In Section 3 we give the proof of
Theorem~\ref{th1}. Section 4 contains some possible directions for
future work.

\section{Background}       

We refer the reader to~\cite{ADFS} for a nice and brief
introduction to Fourier analysis of $\mathbb{Z}_r^n$. In the
following we recall some basic facts and introduce some notations.

Let $r>2$ and $G=\{0,1,\ldots, r-1\}^n=\mathbb{Z}_r^n$. For any $S
\in G$, let $S_i$ denote the $i$th coordinate of $S$. We also think
of $G$ as probability space endowed with the uniform (product)
measure $\mu$.

For any $S\in G$ let $u_S:G\rightarrow \mathbb{C}$ be defined by
$$u_S(T)=\exp\left(\frac{2\pi i\sum_{i=1}^{n} S_iT_i}{r}\right).$$
It is well-known that the set of all functions $u_S$ ($S\in G$)
forms an orthonormal basis for the space of all functions
$f:G\rightarrow \mathbb{C}$. Therefore any such $f$ has a unique
expansion of the form $f=\sum \widehat{f}(S)u_S$, where
$$\widehat{f}(S)=\langle f,u_S\rangle=
\int f(T)\cdot\overline{u_S(T)} \mu(dT).$$
For any function $f:G\rightarrow \mathbb{C}$, define the $p$-norm
of $f$ as
$$\|f\|_p={\left(\int |f(S)|^p \mu(dT)\right)}^{\frac{1}{p}}.$$
From orthogonality it can be easily seen that
$$\|f\|_2^2=\sum_{S\in G}
\widehat{f}(S)^2$$
and
$$\langle
f,g\rangle=\sum\widehat{f}(S)\overline{\widehat{g}(S)}.$$

We use the following notations throughout the paper: For every
complex number $z$, let $d(z,\{0,1\})=\min(|z|,|z-1|)$ denote its
distance from the nearest element in $\{0,1\}$. For any $S\in G$ let
$|S|=|\{i:S_i\neq 0\}|$. $\overline{0}=(0,0,\ldots,0)$, and for each
$1\leq i\leq n$ let $e_i=(0,\ldots,1,\ldots,0)$ be the unit vector
with 1 at $i$th coordinate.  Define $F_S$ as $F_S=
\widehat{f}(S)u_{S}$. Let $f^{>k}=\sum_{|S|>k}F_S$ (similarly
$f^{<k}=\sum_{|S|<k}F_S$) and $f^{=k}=\sum_{|S|=k}F_S$. We
occasionally refer to $f^{=k}$ the $k$-th \emph{level} of Fourier
expansion of $f$. Note that for any function $f$,
$\widehat{f}(\overline{0})$ is the expectation of $f$, and
$\|f^{\geq 1}\|_2^2$ is the variance of $f$.

The following version of Bennett's Inequality which can be easily
obtained from the one stated in~\cite{Bennett} will be used in the
proof of Lemma~\ref{l1} below.

\begin{theorem}[Bennett's Inequality]
\label{Bennett} Let $X_1,\ldots,X_n$ be independent real-valued
random variables with zero mean, and assume that $X_i \le c$ with
probability one. Let
$$\sigma^2=\sum_{i=1}^{n}\Var[X_i].$$
Then for any $t>0$,
$$\Pr[\sum X_i \ge t] \le e^{-\frac{\sigma^2}{c^2} h(\frac{tc}{\sigma^2})},$$
where $h(u)=(1+u)\ln(1+u)-u$ for $u\ge 0$.
\end{theorem}

\section{Main results.}

In~\cite{B1,FKN,KS} results of the following type have been proven:
Let $f$ be a Boolean function on $\mathbb{Z}_2^n$ and $f^{>k}$ is
sufficiently small for some constant $k$, then $f$ is close to being
determined by a few number of coordinates. The following lemma which
is the key lemma in the proof of Theorem~\ref{th1} is a result of
this type for $\mathbb{Z}_r^n$.

\begin{lemma}
\label{l1} Let $f:\mathbb{Z}_r^n \rightarrow \mathbb{C}$ be a
Boolean function such that $\|f^{=1}\|_2^2 \le \frac{1}{r}$ and
$\|f^{>1}\|_2^2\le \epsilon$, where $\epsilon<\frac{1}{10^8r}$ and
$r\ge 20$. Then denoting by $1 \le i_0 \le n$ the index such that
$\sum_{j=1}^{r-1} |\widehat{f}(je_{i_0})|^2$ is maximized, we have

$$\left\|f-\left(\widehat{f}(\overline{0})+\sum_{j=1}^{r-1}
F_{je_{i_0}}\right)\right\|_2^2 < 5\epsilon.$$
\end{lemma}

\begin{remark}
Lemma~\ref{l1} shows that $f$ is close to a function which depends
only on the $i_0$-th coordinate. We do not know if the condition
$\|f^{=1}\|_2^2 \le \frac{1}{r}$ is a weakness of our proof or it is
essential. The condition $\epsilon<\frac{1}{10^8r}$ is not a major
weakness, since for $\epsilon\ge \frac{1}{10^{8}r}$, we have the
trivial bound of $(10^8+1)\epsilon$.
\end{remark}

We postpone the proof of Lemma~\ref{l1} until
Section~\ref{lemmaproof}. We now give the proof of Theorem
\ref{th1}, assuming Lemma~\ref{l1}.

\noindent
\begin{proof}[{\bf Theorem \ref{th1}}]
Let $J$ be an independent set of $G$ such that
$\frac{|J|}{|G|}=\frac{1}{r}(1-\epsilon)$. Let $f$ be the
characteristic function of $J$. Then according to the proof of
Theorem~\ref{t2} (Theorem 1.2 in \cite{ADFS}), we have
$$\|f^{> 1}\|_2^2 = \sum_{|S|>1} {|\widehat{f}(S)|^2} \leq \frac{2\epsilon}{r}.$$
Since
$$\|f^{=1}\|_2^2 \le \|f\|_2^2=\mu(J)\leq \frac{1}{r},$$
by Lemma~\ref{l1}, there exists a function $g:\mathbb{Z}_r^n
\rightarrow \mathbb{C}$ which depends on one coordinate and
$\|f-g\|_2^2\leq \frac{10\epsilon}{r}$. By rounding $g$ to the
nearest of 0 or 1, we get a Boolean function $g_1$ which depends on
one coordinate, and since $f$ is Boolean
$$\|f-g_1\|_2^2\leq 4 \|f-g\|_2^2 \leq \frac{40\epsilon}{r}.$$
\end{proof}

\subsection{Proof of Lemma~\ref{l1} \label{lemmaproof}}

The proof of Lemma~\ref{l1} shares similar ideas with the proof of
Theorem 8 in~\cite{KS}. However dealing with (complex) Fourier
expansions on $\mathbb{Z}_r^n$ instead of (real) generalized Walsh
expansions on $\mathbb{Z}_2^n$ required new ideas.

For any function $f$, denote $\widehat{f}(S)u_S$ by $F_S$ to make
the notations easier. For $1 \le i \le n$, let
$g_i=\sum_{j=1}^{r-1}F_{je_i}$, and define
$g_0=\widehat{f}(\overline{0})$. For $0 \le i \le n$ let
$a_i=\|g_i\|_2$. Without loss of generality assume that $a_1 \ge a_2
\ge \ldots \ge a_n$. To obtain
$$\|f-(g_0+g_1)\|_2^2=\sum_{i=2}^{n} a_i^2 + \|f^{>1}\|_2^2 \leq
5\epsilon,$$
we will first show that $a_2$ is small (Claim~\ref{claim1}). This
would allow us to apply a concentration theorem and conclude that
$\sum_{i=2}^{n} a_i^2$ is very small (Claim~\ref{claim2}).

First note that
$$\|f^{=1}\|_2^2=\sum_{i=1}^n a_i^2 \le \frac{1}{r},$$
which implies that $a_2^2 \le \frac{1}{2r}$. Now since $\|g_2\|_2^2
\le \frac{1}{2r}$, for every $0\le x_2 \le r-1$,
\begin{equation}
\label{g_2} |g_2(x_2)|\le \sqrt{1/2}.
\end{equation}
\begin{claim}
\label{claim1}$a_2^2 < 2000\epsilon$.
\end{claim}
\begin{proof}
Consider an arbitrary assignment $\delta_1,\delta_3,\ldots,\delta_n$
to $x_1,x_3,\ldots,x_n$, and let
$$l=\widehat{f}({\overline{0}})+g_1(\delta_1)+\sum_{i=3}^n g_i(\delta_i).$$
Since for every $0 \le x_2 \le r-1$,
$$d(l,\{0,1\}) \le |g_2(x_2)|+d(l+g_2(x_2),\{0,1\}),$$
we have
$$\|d(l,\{0,1\})\|_2^2 \le
2(\|g_2\|_2^2+\|d(l+g_2,\{0,1\})\|_2^2),$$
or equivalently
\begin{equation}
\label{bound1} d(l,\{0,1\})^2 \le 2(a_2^2+\|d(l+g_2,\{0,1\})\|_2^2).
\end{equation}
Note that
$$\|d(f^{\le 1},\{0,1\})\|_2^2 \le 2(\|d(f,\{0,1\})\|_2^2+\|f^{>1}\|_2^2) \le 2\epsilon.$$
Therefore we can find an assignment
$\delta_1,\delta_3,\ldots,\delta_n$ such that

\begin{equation}
\label{boundlg} \|d(l+g_2,\{0,1\})\|_2^2\leq 2\epsilon.
\end{equation}

By~(\ref{bound1}) for any such assignment, we have
$d(l,\{0,1\})^2\le\frac{1}{r}+4\epsilon\le 1/16$, which implies
either $|l|\le \frac{1}{4}$ or $|l-1|\le \frac{1}{4}$.

Define
$\lambda=\frac{1-\sqrt{\frac{1}{2}}-\frac{1}{4}}{\sqrt{\frac{1}{2}}+\frac{1}{4}}$.
Now~(\ref{g_2}) implies that for any $0\leq x_2\leq r-1$,
\begin{itemize}
\item[{\bf Case 1:}] If $|l|<\frac{1}{4}$, then
$|(l+g_2(x_2))-1|\ge \lambda |l+g_2(x_2)|$.

\item[{\bf Case 2:}] If $|l-1|<\frac{1}{4}$, then
 $|l+g_2(x_2)|\ge
\lambda |(l+g_2(x_2))-1|$.
\end{itemize}

Let $A=\{x_2\in \mathbb{Z}_r:|l+g_2(x_2)|\leq|l+g_2(x_2)-1|\}$ and
denote its complement by $\overline{A}$. Representing
$\|d(l+g_2,\{0,1\})\|_2^2$ as a sum of two integrals over $A$ and
$\overline{A}$, and using~(\ref{g_2}), in Cases~1 and~2 one can show
that
$$\|d(l+g_2,\{0,1\})\|_2^2 \ge
\lambda^2\|g_2\|_2^2> \frac{a_2^2}{1000}.$$
Note that the assumption $a_2^2 \ge 2000\epsilon$ will imply
$\|d(l+g_2,\{0,1\})\|_2^2>2\epsilon$ which
contradicts~(\ref{boundlg}). Thus $a_2^2 < 2000\epsilon$.
\end{proof}

\begin{claim}
\label{claim2} $\sum_{i=2}^n a_i^2 \le 4 \epsilon$.
\end{claim}
\begin{proof}
Let $2 \le m \le n$ be the minimum index which satisfies
\begin{equation}
\label{define_m} \sum_{i=m}^n a_i^2 \le 10^4 \epsilon.
\end{equation}

Denote $I=\{m, \ldots, n\}$, and for every $y \in
\mathbb{Z}_r^{m-1}$ let $f^*_{I[y]}$ be a function of
$\mathbb{Z}_r^{n-m+1}$ (with uniform measure $\mu'$) defined as
$$f^*_{I[y]}(x)=f^{\le 1}(y \cup x).$$
Obviously
$$\int \|d(f^*_{I[y]}(x),\{0,1\})\|_2^2 \mu'(dy) = \|d(f^{\le 1},\{0,1\})\|_2^2 \le 2\epsilon.$$
Hence for some $y$, $\|d(f^*_{I[y]}(x),\{0,1\})\|_2^2 \le
2\epsilon$. Let $b=\widehat{f}(\overline
0)+\sum_{i=1}^{m-1}g_i(y_i)$. Then
$$f^*_{I[y]}(x)=b+\sum_{i=m}^n g_i(x_i).$$
Applying Lemma~\ref{concent} below to $f^*_{I[y]}$ for
$\epsilon'=2\epsilon$ shows that $\sum_{i=m}^n a_i^2 \le 4
\epsilon$. This will imply that $m=2$, as $a_2^2<2000 \epsilon$ and
$m$ was the minimum index satisfying~(\ref{define_m}), which
completes the proof.
\end{proof}

\begin{lemma}
\label{concent} Let $f:\mathbb{Z}_r^n \rightarrow \mathbb{C}$ be a
function satisfying $f^{>1} \equiv 0$. Let $ \| d(f,\{0,1\})
\|_2^2 \le \epsilon'$, and suppose that $\| f^{=1}\|_2^2 <
10^4\epsilon'$ and $\epsilon'<\frac{2}{10^8r}$. Then we have
$$\|f^{=1}\|_2^2<2\epsilon'.$$
\end{lemma}
\begin{proof}
Suppose that $f=b+\sum_{i=1}^{n} g_i$, where
$b=\widehat{f}(\overline{0})$ and $g_i=\sum_{j=1}^{r-1} F_{j e_i}$.
We have
$$\|d(b,\{0,1\}) \|_2^2 \le 2(\|d(f,\{0,1\})\|_2^2 + \|f-b\|_2^2) \le
20002\epsilon'.$$
Without loss of generality assume that
$d(b,1)\le\sqrt{20002\epsilon'}$ which implies that
\begin{equation}
\label{boundB} \Re(b)>2/3.
\end{equation}
We have
$$\|f-1\|_2^2-\|d(f,\{0,1\})\|_2^2 = \int (|f-1|^2-|f|^2)\zeta dx,$$
where
$$\zeta(x)=\left\{
\begin{array}{cr} 1 & \Re(f(x))<\frac{1}{2} \\
0 & {\rm otherwise}
\end{array}\right.$$
So
\begin{equation}
\label{difference} \|f-1\|_2^2-\|d(f,\{0,1\})\|_2^2 = \int
(1-2\Re(f))\zeta dx.
\end{equation}

The next step is to show that (\ref{difference}) is less than
$\epsilon'$. Note that $\Re(f)=\Re(b)+\sum_{i=1}^n \Re(g_i)$, and
$\int \Re(g_i)=0$. Moreover
$$\int \Re(g_i)^2 = \|\Re(g_i)\|_2^2 \le \|g_i\|_2^2.$$
So
$$ \|\Re(g_i)\|_2^2 \le \sum \|g_i\|_2^2 \le 10^4 \epsilon'$$
which follows that for every $x$,
$$|\Re(g_i(x))| \le \sqrt{10^4 r \epsilon'} \le
\sqrt{2}\times10^{-2} \doteq c.$$
Applying Theorem~\ref{Bennett} with $X_i=-\Re(g_i)$, we get
\begin{equation}
\Pr[\sum \Re(g_i) \le -t] \le e^{\frac{-10^4 \epsilon'}{c^2}
h(10^{-4}tc/\epsilon')}.
\end{equation}
Note that $h(u) \ge u \ln(\frac{u}{e})$, for $u \ge e$; which
implies that for $t \ge \frac{1}{6} \ge 10^4 e \epsilon'/c$,
\begin{equation}
\label{concenteration} \Pr[\sum \Re(g_i) \le -t] \le e^{-\frac{t}{c}
\ln(10^{-4}tc/ e \epsilon')}.
\end{equation}
Now
$$(\ref{difference}) =
\int_{t=0}^{\infty} \Pr[1-2\Re(f)>t] =\int_{t=0}^{\infty}
\Pr\left[\Re(b)+\sum \Re(g_i)<\frac{1-t}{2}\right].$$
Substituting~(\ref{boundB}) we get
$$(\ref{difference}) \le  \int_{t=0}^{\infty} \Pr\left[\sum
\Re(g_i)<\frac{1-t}{2}-\frac{2}{3}\right] \le 2\int_{t=1/6}^{\infty}
\Pr\left[\sum \Re(g_i)<-t\right].$$
Now by (\ref{concenteration})
\begin{eqnarray}
\label{finaldiff} (\ref{difference}) &\le& 2\int_{t=1/6}^{\infty}
e^{\frac{t}{c} \ln(10^4 e \epsilon' /t c)}\le 2\int_{t=1/6}^{\infty}
\left(\frac{1-\ln(10^4e\epsilon'/t c)}{c}\right) e^{\frac{t}{c}
\ln(10^4 e \epsilon'/t c)} = \nonumber
\\ &= & 2e^{\frac{1}{6c} \ln(6\times10^4 e\epsilon'/c)} <\epsilon',
\end{eqnarray}
because $\epsilon'\le 10^{-8}$. Finally by~(\ref{finaldiff})
$$\|f^{=1}\|_2^2 \le \|f-1\|_2^2 \le \|d(f,\{0,1\})\|_2^2 +
\epsilon' \le 2\epsilon'.$$
\end{proof}

\section{Future Directions}

Lemma~\ref{l1} asserts that when most of the 2-norm weight of the
Fourier expansion of a Boolean function on $\mathbb{Z}_r^n$ is
concentrated on the first two levels, then the function can be
approximated by a Boolean function that depends on only one
coordinate. One possible generalization of this lemma would be to
show that a Boolean function on $\mathbb{Z}_r^n$ whose Fourier
expansion is concentrated on the first $l$ levels for some constant
$l$ can be approximated by a Boolean function that depends on $k(l)$
coordinates, for some function $k(l)$. Analogues of this for
$\mathbb{Z}_2^{n}$ have been proven in~\cite{B1} and~\cite{KS}.

Consider a graph $G$ whose vertices are the elements of the
symmetric group $S_n$ and two vertices $\pi$ and $\pi'$ are adjacent
if $\pi(i)\neq \pi'(i)$ for every $1 \le i \le n$. For every $1 \le
i,j \le n$ the set $S_{ij}$ of the vertices $\pi$ satisfying
$\pi(i)=j$ forms an independent set of size $(n-1)!$. Recently
Cameron and Ku~\cite{Cameron} have proved that these sets are the
only maximum independent sets of this graph. Similar results have
been proven for generalizations of this graph in~\cite{LM}. Cameron
and Ku made the following conjecture:

\begin{conj}\label{Cameron}
{\bf \cite{Cameron}} There is a constant $c$ such that every
independent set of size at least $c(n-1)!$ is a subset of an
independent set of size $(n-1)!$.
\end{conj}

One might notice the similarity of Conjecture~\ref{Cameron} and
Corollary~\ref{cor1} for $r=n$. Despite this similarity we are not
aware of any possible way to apply the techniques used in this paper
to the problem.  Since $S_n$ is not Abelian, the methods of the
present paper (and all the papers mentioned in Section~\ref{intro})
fail to apply directly to this problem. So an answer to
Conjecture~\ref{Cameron} or its analogues for the graphs studied
in~\cite{LM} (which do not even have a group structure) might lead
to new techniques.

\bibliographystyle{siam}
\bibliography{independent}
\end{document}